\documentclass[11pt, oneside]{article}   	% use "amsart" instead of "article" for AMSLaTeX format
\usepackage{geometry}                		% See geometry.pdf to learn the layout options. There are lots.
\usepackage{graphicx}				% Use pdf, png, jpg, or eps§ with pdflatex; use eps in DVI mode
\usepackage{amsmath,amssymb,amsfonts,graphicx}
\usepackage{multirow,color}

\newtheorem{lemma}{Lemma}[section]

\newtheorem{theorem}[lemma]{Theorem}

\newcommand\blfootnote[1]{%
  \begingroup
  \renewcommand\thefootnote{}\footnote{#1}%
  \addtocounter{footnote}{-1}%
  \endgroup
}

\usepackage{lipsum}
\usepackage{amsfonts}
\usepackage{graphicx}
\usepackage{epstopdf}

\def\e{\varepsilon}

\title{Asymptotic behaviour of the capacity\\ in two-dimensional heterogeneous media}
\author{A. Braides and G.C. Brusca
\\ 
SISSA, via Bonomea 265, Trieste, Italy}
\date{}							% Activate to display a given date or no date

\begin{document}
\maketitle

%\hfill{\em Dedicated to Irene Fonseca on the occasion of her 65th birthday}

\begin{abstract}
We describe the asymptotic behaviour of the minimal inhomogeneous two-capacity of small sets in the plane with respect to a fixed open set $\Omega$. This problem is governed by two small parameters: $\e$, the size of the inclusion (which is not restrictive to assume to be a ball), and $\delta$, the period of the inhomogeneity modelled by oscillating coefficients. We show that this capacity behaves as $C|\log\e|^{-1}$. The coefficient $C$ is explicitly computed from the minimum of the oscillating coefficient and the determinant of the corresponding homogenized matrix, through a harmonic mean with a proportion depending on the asymptotic behaviour of $|\log\delta|/|\log\e|$.
\smallskip

    
{\bf AMS Class:} 49J45, 35B27, 31A15, 74Q05
\end{abstract}

\blfootnote{Preprint SISSA 10/2022/MATE}
\ \bigskip

\hfill{\em Dedicated to Irene Fonseca on the occasion of her 65th birthday}

\section{Introduction}
%\subsection{}
Scaling-invariant variational problems present challenging technical issues and intriguing properties due to the concentration of minimal energy configurations on several scales. The paradigmatic example highlighting the effects of scaling invariance can be considered that of Ginzburg-Landau functionals, whose minimizers may generate vortices triggered by topological incompatibilities (see e.g.~\cite{BBH}). A recent paper \cite{ABCDP} investigates the minimal energy of Ginzburg-Landau vortices in heterogeneous media in dimension two, showing the different type of interactions between the Ginzburg-Landau energy and heterogeneities at different scales. An even simpler scaling-invariant functional is the Dirichlet integral in dimension two, which entails singular properties of the two-capacity (more in general, this holds for the $d$-capacity in dimension $d$), and for which in this paper we examine a problem analogous to that of \cite{ABCDP}, highlighting the corresponding analogies and differences.

\smallskip
We consider the problem of the minimal capacity of a small set in a highly heterogeneous media at the critical scaling. We restrict to the planar case in which we may provide a clearer explanation of the asymptotic analysis confining to quadratic functionals.
It is not restrictive to suppose that the inclusion is a ball, for which the problem can be stated as the study of the asymptotic behaviour of minima
\begin{equation}\label{emmedeltae}
m_{\e,\delta}=\min\Bigl\{\int_{\Omega} a\Bigl({x\over\delta}\Bigr)|\nabla u|^2\,dx: u\in H^1_0(\Omega), u=1 \hbox{ on } B_\e(z), z\in\Omega\Bigr\},
\end{equation} 
with $\Omega$ a fixed bounded open set of $\mathbb R^2$. 
Here, $a$ is a $1$-periodic function representing the geometry of heterogeneities, $\delta$ is the scale of oscillations of the heterogeneities and $B_\e(x)$ is the small inclusion, a ball of centre $x$ and radius $\e$. Since the asymptotic result is independent of the particular $\Omega$ we omit the dependence on the set in the notation in \eqref{emmedeltae}. 

Note that if $a$ in \eqref{emmedeltae} is a constant $c$ the value of $m_\e=m_{\e,\delta}$ can be easily computed recalling that for $R>r$ we have
\begin{equation}\label{eqRr}
\min\Bigl\{\int_{B_R}|\nabla u|^2\,dx: u\in H^1_0(B_R), u=1 \hbox{ in } B_r\Bigr\}= 2\pi{1\over \log(R/r)},
\end{equation}
where $B_R=B_R(0)$ and $B_r=B_r(0)$.
This computation is straightforward since the minimizer is an explicit logarithmic function, and implies that 
\begin{equation}\label{casoii-7}
\lim_{\e\to 0}|\log\e| m_{\e}=
2\pi c\,.
\end{equation}

 We now turn to a non-constant coefficient $a$. We set $\alpha=\hbox{\rm ess-}\inf a>0$ and let $A_{\rm hom}$ denote the constant-coefficient matrix of the corresponding homogenized quadratic form.
%
%Note that if $\e$ is fixed then $m_{\e,\delta}$ converges to 
%$$
%m_{\e}=\min\Bigl\{\int_{\Omega} \sqrt{\alpha\beta}|\nabla u|^2\,dy: u\in  W^{1,d}_0(\Omega), u=1 \hbox{ on } B_\e(x), x\in\Omega\Bigr\}.
%$$
%The behaviour of $m_\e$ is logarithmic, in the sense that 
%\begin{equation}\label{limbeh}
%\lim_{\e\to 0}|\log\e| m_{\e}=2\pi\sqrt{\alpha\beta}.
%\end{equation}
We consider $\delta=\delta_\e$ and study the asymptotic behaviour of $m_{\e,\delta}$ as $\e,\delta\to 0$.
We will prove that the limit analysis is determined by the parameter $\lambda\in[0,1]$ defined by
\begin{equation}\label{lambda}
\lambda=\lim_{\e\to0}{|\log \delta|\over|\log\e|}\wedge 1
\end{equation}
(a limit which we may suppose exists up to subsequences), and that
the behaviour of $m_{\e,\delta}$ is logarithmic; more precisely, that we have 
\begin{equation}\label{casoii}
\lim_{\e\to 0}|\log\e| m_{\e,\delta}=%\begin{cases} \alpha &\hbox{ if } \lambda=0\\
2\pi {\alpha\sqrt{{\rm det}A_{\rm hom}}\over \lambda \alpha+(1-\lambda)\sqrt{{\rm det}A_{\rm hom}}} . 
%&\hbox{ if } 0<\lambda<1\\C_{\rm hom} &\hbox{ if } \lambda=1\end{cases}
\end{equation}
Recalling  \eqref{casoii-7}, we remark that \eqref{casoii} mixes the two cases $c=\alpha$ (minimal value) and $c=\sqrt{{\rm det}A_{\rm hom}}$ (homogenized value) through their {\em harmonic mean}.

We briefly explain formula \eqref{casoii} by describing an approximate minimizer $u^\e$ for $m_{\e,\delta}$ for $\e$ small.
 For simplicity, we suppose $0<\lambda<1$ and $a$ continuous.
We choose the center of the ball $z=z^\e$ to be a minimum point for $a(x/\delta)$, and we fix any $\lambda'<\lambda<\lambda''$; note that $\e^{\lambda''}=:\delta''<\!<\delta<\!< \delta':=\e^{\lambda'}$. We fix a constant $\overline u$ (to be optimized later) and take $u^\e$ on $B''=B_{\delta''}(z^\e)$ defined as the minimizer of 
$$
\min\Bigl\{\int_{B''} \alpha|\nabla u|^2\,dx: u= \overline u\hbox{ on }\partial B'', u=1 \hbox{ in } B_\e(z^\e)\Bigr\}= 2\pi \alpha{1\over |\log\e|(1-\lambda'')} (1-\overline u)^2
$$
by \eqref{eqRr}. Note that in $B''$ we have $a(x/\delta)\sim\alpha$.
We may also take $u^\e=\overline u$ on $B'\setminus B''$, where $B'=B_{\delta'}(x^\e)$. 
Indeed, letting $\lambda''-\lambda'\to 0$ we see that the contribution of the Dirichlet integral on $B'\setminus B''$ is negligible. 
We may suppose that $B_1(z^\e)\subset \Omega$, so that we can take $u^\e$ on $B_1(z^\e)\setminus \overline{B'}$ as the function minimizing 
$$
\min\Bigl\{\int_{B_1(z^\e)} a\Bigl({x\over \delta}\Bigr)|\nabla u|^2\,dx: u\in H^1_0(B_1(x^\e)), u=\overline u\hbox{ in } B' \Bigr\}.
$$
By the scale-invariance properties of the Dirichlet integral in two dimension, we can scale the problem to
$$
\min\Bigl\{\int_{B_{1\over\delta'}({z^\e\over\delta'})} a\Bigl({x\over \delta/\delta'}\Bigr)|\nabla u|^2\,dx: u\in H^1_0\Bigl(B_{1\over\delta'}\Bigl({z^\e\over\delta'}\Bigr)\Bigr), u=\overline u\hbox{ in } B_1\Bigl({z^\e\over\delta'}\Bigr) \Bigr\}.
$$
Since $\delta/\delta'<<1$, this problem can be approximated by its homogenized one
$$
\min\biggl\{\int_{B_{1\over\delta'}({z^\e\over\delta'})} \langle A_{\rm hom}\nabla u, \nabla u\rangle\,dx: u\in H^1_0\Bigl(B_{1\over\delta'}\Bigl({z^\e\over\delta'}\Bigr)\Bigr), u=\overline u\hbox{ in } B_1\Bigl({z^\e\over\delta'}\Bigr) \biggr\}\sim {2\pi \sqrt{\hbox{\rm det}A_{\rm hom}}\over |\log\e|\lambda'} \overline u^2
$$
by \eqref{eqRr} and a change of variables. Summing up these two values and letting $\lambda',\lambda''\to \lambda$ we have the estimate 
$$
\limsup_{\e\to0} m_{\e,\delta}|\log\e|\le 2\pi \Bigl(\sqrt{\hbox{\rm det}A_{\rm hom}} {1\over\lambda}\overline u^2+\alpha{1\over1-\lambda}(1-\overline u)^2\Bigr)
$$
Taking $\overline u$
%={\alpha\lambda\over \sqrt{\alpha\beta}(1-\lambda) +\alpha\lambda}$ 
minimizing the last expression we obtain \eqref{casoii}. Note that in the argument above some care must be used in the passage to the limit for the homogenization since it is performed on varying domains. The optimality of this construction is rather technical and involves a lemma that allows to decompose minimizations for $m_{\e,\delta}$ into separate minimization on annuli of given ratio between inner and outer radii.

%We treat different regimes.
%
%(i) The case $\delta<\!<\e$. In this case we prove that the asymptotic behaviour is the same as in the case of in which we regard $\e$ fixed; more precisely, we prove that
%\begin{equation}\label{casoi}
%\lim_{\e\to 0}|\log\e| m_{\e,\delta}=C_{\rm hom},.
%\end{equation}
%
%
%(ii) The case $\delta\ge C\e$ for some $C>0$. In this case we can assume that there exists the limit
%$$
%\lambda=\lim_{\e\to0}{|\log \delta|\over|\log\e|} \in [0,1].
%$$
%We then prove that 
%\begin{equation}\label{casoii}
%\lim_{\e\to 0}|\log\e| m_{\e,\delta}=%\begin{cases} \alpha &\hbox{ if } \lambda=0\\
%2\pi {\alpha\sqrt{\alpha\beta}\over \lambda \alpha+(1-\lambda)\sqrt{\alpha\beta}} ; 
%%&\hbox{ if } 0<\lambda<1\\C_{\rm hom} &\hbox{ if } \lambda=1\end{cases}
%\end{equation}
%that is, the harmonic mean of $\alpha$ and $\sqrt{\alpha\beta}$ with proportion $\lambda$.

\smallskip
We can compare the result above with the one in \cite{ABCDP}, regarding the asymptotic analysis of the heterogeneous Ginzburg-Landau energy
\begin{equation}\label{GLe}
\int_{\Omega} a\Bigl({x\over\delta}\Bigr)|\nabla u|^2\,dx+{1\over\e^2}\int_\Omega (|u|^2-1)^2\,dx
\end{equation}
at the vortex scaling, where $u:\Omega\to \mathbb R^2$. Such an analysis can be formalized as the study of the asymptotic behaviour of the {\em hard-core energy}
\begin{equation}\label{emmedeltaegl}
m^{\rm hc}_{\e,\delta}=\min\Bigl\{\int_{\Omega} a\Bigl({x\over\delta}\Bigr)|\nabla u|^2\,dx: u\in H^1(\Omega; S^1), {\rm deg}\, u=1 \hbox{ on } \partial B_\e(z), z\in\Omega\Bigr\},
\end{equation} 
leading to an analogy with problem \eqref{emmedeltae}, with a condition on the degree of the test functions
in the place of fixing boundary values.
The behaviour is determined by the same parameter $\lambda$ in \eqref{lambda}, but in this case we have
\begin{equation}\label{casoiigl}
\lim_{\e\to 0}|\log\e| m^{\rm hc}_{\e,\delta}=
2\pi \Bigl(\lambda \sqrt{\hbox{\rm det}A_{\rm hom}}+(1-\lambda)\alpha\Bigr); 
\end{equation}
that is, $2\pi$ times the {\em arithmetic mean} of $\sqrt{\hbox{\rm det}A_{\rm hom}}$ and $\alpha$ with proportion $\lambda$. Again the minimization can be split in two regions, between $B_\e(z^\e)$ and $B''$ and between $B'$ and $B_1(z^\e)$, with a condition on the degree. Minimization in the first region considers for $a$ only the value $\alpha$, and minimization in the second region can be approximated by substituting the homogenized quadratic form $\langle A_{\rm hom}\xi,\xi\rangle$ to the oscillating quadratic form  $a(x/\delta)|\xi|^2$. Both problems behave logarithmically in complete analogy with \eqref{eqRr}. Contrary to the capacitary problem, here there is no further optimization in the boundary condition, and the two minima can be simply added, to obtain \eqref{casoiigl}.

Note that in both problems we have `well-separated' extreme regimes when
$$
\lim_{\e\to 0}|\log\e| m_{\e,\delta}=\begin{cases} 2\pi \alpha &\hbox{ if }|\log\delta|<\!<|\log\e|\\
%{C_{\rm hom}\alpha\over \lambda \alpha+(1-\lambda)C_{\rm hom}}  &\hbox{ if } 0<\lambda<1
2\pi \sqrt{{\rm det}A_{\hom}} &\hbox{ if }{|\log\delta|\over|\log\e|}\ge 1+o(1).\end{cases}
$$
In these regimes we have a ``separation of scales'': we can formally optimize the capacity (or compute the vortex energy) at fixed $\delta$ or apply the homogenization procedure first, respectively, while in all other regimes the two effects both appear. This is in contrast with problems in the gradient theory of phase transitions, with energies formally similar to \eqref{GLe} but with scalar $u$, where separation of scales occurs for all regimes of $\e$ and $\delta$ except when $\e\sim\delta$ (see \cite{ABCP,FC,Morfe}).

\section{Asymptotic analysis}
We assume $\delta=\delta_{\e}$, meaning that $\delta$ is a function of $\e$ defined on $(0,1)$, which vanishes as $\e$ tends to $0$. Up to subsequences, we may suppose that there exists the limit $\lim_{\e\to0} \frac{|\log\delta|}{|\log\e|}$, and we define
\begin{equation}\label{lambda-2}
\lambda=\min\Bigl\{1,\lim_{\e\to0} \frac{|\log\delta|}{|\log\e|}\Bigr\}\in [0,1]\,
\end{equation}
as in \eqref{lambda}
Note that if $\delta\le C\e$ then $\lambda=1$, while if $\delta=\e^\eta$ with $\eta\in (0,1)$ then $\lambda=\eta$. The cases $\lambda=0$ and $\lambda=1$ include also very slowly converging $\delta$, e.g.~$\delta=|\log\e|^{-1}$, and almost linear $\delta$, e.g.~$\delta=\e|\log\e|$, respectively. 

We will simplify our exposition by choosing a special form of $a$; namely, we let $a:\,\mathbb{R}^2 \rightarrow \mathbb{R}$ be a $1$-periodic checkerboard function with values $\alpha, \beta$, such that $0<\alpha<\beta$, defined on $[0,1)^2$ as
\begin{equation}\label{fun-a}
a(x):=
\begin{cases}
\alpha & \hbox{ if } x\in \left[0,\frac{1}{2}\right)^2 \cup \left[\frac{1}{2},1\right)^2 \\
\beta & \hbox{ if } x\in \left[\frac{1}{2},1\right)\times \left[0,\frac{1}{2}\right) \cup \left[0, \frac{1}{2}\right)\times \left[\frac{1}{2},1\right).
\end{cases}
\end{equation}
The analysis of this particular inhomogeneity will contain all the main features of the general case, with some simplifications in the proofs. 
We will prove the following result.

\begin{theorem}[asymptotic behaviour of the capacity in a heterogeneous medium]\label{main} Let $\Omega$ be a bounded open set in $\mathbb R^2$, and, for $\e, \delta>0$, let $m_{\e,\delta}$ be defined by
\begin{equation}\label{emmedeltae-0}
m_{\e,\delta}=\min\Bigl\{\int_{\Omega} a\Bigl({y\over\delta}\Bigr)|\nabla u|^2\,dy: u\in H^1_0(\Omega), u=1 \hbox{ on } B_\e(x), x\in\Omega\Bigr\}.
\end{equation} 
If $\delta=\delta_\e$ is infinitesimal as $\e\to0$, and $\lambda$ is defined as in \eqref{lambda-2}, then we have 
\begin{equation}\label{casoii-0}
\lim_{\e\to 0}|\log\e| m_{\e,\delta}=%\begin{cases} \alpha &\hbox{ if } \lambda=0\\
2\pi {\alpha\sqrt{\alpha\beta}\over \lambda \alpha+(1-\lambda)\sqrt{\alpha\beta}} . 
%&\hbox{ if } 0<\lambda<1\\C_{\rm hom} &\hbox{ if } \lambda=1\end{cases}
\end{equation}
\end{theorem} 

The proof of this theorem will be split between a lower and an upper bound in the following sections.

The appearance of the term $\sqrt{\alpha\beta}$ is a consequence of a classical homogenization result (see \cite{JKO,BDF}). For convenience in the proofs of the following sections, we state this result below with some translations $\tau_\eta$. To formally derive this modified statement from the usual one without translations it suffices first to note that by the periodicity of $a$ we can suppose $\|\tau_\eta\|_\infty\le 1$, so that if $u_\eta\to u$ in $L^2(\mathbb R^2)$ then, defining $\overline u_\eta(x)= u_\eta(x-\eta\tau_\eta)$, we still have $\overline u_\eta\to u$ in $L^2(\mathbb R^2)$. This change of variables allows to deduce the modified statement in the following form.

\begin{theorem}[homogenization of the square checkerboard]\label{homab} Let $A$ be a bounded open subset of $\mathbb R^2$ with Lipschitz boundary. For $\eta>0$ let $\tau_\eta\in\mathbb R^2$ let 
$$
F_\eta(u)= \int_A a\Bigl({x\over\eta}+\tau_\eta\Bigr)|\nabla u|^2\,dx\qquad u\in H^1(A).
$$
Then $F_\eta$ $\Gamma$-converge to $\sqrt{\alpha\beta}\int_A|\nabla u|^2\,dx$ with respect to the $L^2(\Omega)$-convergence. In particular, for all $\phi\in H^1(A)$ we have 
$$
\lim_{\eta\to 0} \min\{F_\eta(u) : u=\phi\hbox{ on }\partial A\}=\sqrt{\alpha\beta}\min\Bigl\{\int_A |\nabla u|^2\,dx: u=\phi\hbox{ on }\partial A\Bigr\}.
$$
Furthermore this convergence is uniform in $\tau_\eta$.
\end{theorem}

A fundamental tool in the proof of Theorem \ref{main} is a variant elaborated in \cite{AnsBra} (see also \cite{Sigalotti}) of a method by De Giorgi to vary boundary conditions on converging sequences (see \cite{GCB,DM}). In our context it can be stated as follows.

\begin{lemma}[reduction to functions constant on a suitable circumference]\label{modlemma1}
Let $\gamma: \mathbb{R}^2\rightarrow \mathbb{R}$ be a measurable function and let $\alpha,\beta$ be real constants such that $0<\alpha\leq\gamma(x)\leq\beta$ for every $x\in \mathbb{R}^2$. Let $z\in\mathbb{R}^2$, $R>0$ and define
\[
F(u,A):=\,\int_{A} \gamma(x)|\nabla u(x)|^2\,dx
\]
for every $u\in H^1(B_R(z))$ and  measurable subset $A\subseteq B_R(z)$.

We fix $\eta\in (0,1)$ and set $S:=$max$\,\{s\in \mathbb{N} : \eta2^s \leq R\}$. We assume that $S\geq3$.
Let $N$ be a natural number such that $2\leq N<S$ and let $r$ be a positive real number such that $r \leq \eta2^{S-N}$.

Then there exists a function $v$ with the following properties:

{\rm(i)} $v\in H^1(B_R(z)\setminus \overline B_r(z))$, 

{\rm(ii)} there exists $j\in\{1,...,N-1\}$ such that 
$$
v=u \hbox{ on } \bigl(B_{\eta2^{S-j-1}}(z)\setminus \overline B_R(z)\bigr)\cup\bigl( B_R(z)\setminus \overline B_{\eta2^{S-j+1}}(z)\bigr),
$$

{\rm(iii)} for the same $j$ the function $v$  is constant on $\partial B_{\eta2^{S-j}}(z)$,

{\rm(iv)} there exists a positive constant $C$ depending only on $\alpha$ and $\beta$ such that

$$F(v,B_R(z)\setminus \overline B_r(z) )\leq \Bigl(1+\frac{C}{N-1}\Bigr)\, F(u,B_R(z)\setminus \overline B_r(z)).$$

\end{lemma}

This lemma states that, up to a change in the energy, which by (iv) is small if $N$ is large,  by (iii) we can suppose that functions $u$ have a constant value on the boundary of some ball, which is of relative radius close to $R$ if $N$ is small compared to $S$. In the following this will be applied with fixed $N$ and diverging $S$. The proof of this lemma is obtained by fixing $j$, modifying $u$ only on the annulus $B_{\eta2^{S-j+1}}(z)\setminus B_{\eta2^{S-j-1}}(z)$ in such a way that (iii) holds. This is done by a cut-off joining $u$ and its average on an annulus, and estimating the change of energy using a Poincar\'e inequality. By the scaling properties of the Poincar\'e inequality this estimate can be shown to involve a constant $C$ depending only on $\alpha$ and $\beta$ since the sets in which it is applied are all homothetic. Eventually, De Giorgi's trick consists in choosing an optimal $j$ in (ii), which gives estimate (iv). We refer to \cite{Brusca} for an explicit construction.

\subsection{Lower bound}
Let $R_\Omega$ denote the maximum between the diameter of $\Omega$ and $1$ (this condition is just to have $\log R_\Omega$ non negative in the following), so that for all $z\in \Omega$ the extension of any $u\in H^1_0(\Omega)$ belongs to $H^1_0(B_{R_\Omega}(z))$.
By the estimate $a(x)\ge \alpha$ for all $x$ we have 
\begin{eqnarray*}
m_{\e,\delta}&\ge& \alpha\min\Bigl\{\int_{B_{R_\Omega}(z)} |\nabla u|^2\,dx: u\in H^1_0(B_{R_\Omega}(z)), u=1 \hbox{ on } B_\e(z), z\in\Omega\Bigr\}\\
&=& \alpha\min\Bigl\{\int_{B_{R_\Omega}} |\nabla u|^2\,dx: u\in H^1_0(B_{R_\Omega}), u=1 \hbox{ on } B_\e\Bigr\}.
\end{eqnarray*}
Using \eqref{eqRr} we then get 
$$
\liminf_{\e\to0}|\log\e|m_{\e,\delta}\ge 2\pi\alpha,
$$
which is the desired lower bound for $\lambda=0$.  
\medskip

If $\lambda\in (0,1]$ we choose $0<\lambda_1<\lambda$ such that 
\begin{equation}\label{eqz}
\lim\limits_{\e\to0}{\delta\over \e^{\lambda_1}}=0\,.
\end{equation}
We consider $z\in\Omega$ and $u\in H^1_0(\Omega)$ with $u=1$ on $B_\e(z)$, extended to $0$ on $\mathbb R^2\setminus\Omega$. 
 We now construct a family of concentric annuli centered at $z$, to each of which we will apply Lemma \ref{modlemma1}. We first set
\begin{equation}\label{T}
T=\max \{n\in \mathbb{N} : \e^{\lambda_1}2^n \leq R_\Omega\}.
\end{equation}
Note that the set of such $n$ is not empty as soon as $\e$ is small enough. Moreover, we can explicitly write 
$$
T=\Bigl\lfloor\frac{\lambda_1 |\log \e|+\log R_\Omega}{\log 2}\Bigr\rfloor.
$$
In the constructions below we suppose $\e$ small enough so that $T\geq4$. 
We fix natural s $M\in(2,T)$ and $N\in\mathbb N\cap(0,M)$, and apply Lemma \ref{modlemma1} $\lfloor T/M\rfloor$ times, with 
$$
\gamma(x)=a\left(\frac{x}{\delta}\right),\, \eta=\e^{\lambda_1},\, R=\e^{\lambda_1}2^{kM},\, r=\e^{\lambda_1}2^{(k-1)M}\hbox{and } S=kM,
$$
obtaining functions $v^k\in H^1(B_{\e^{\lambda_1}2^{kM}}(z)\setminus \overline B_{\e^{\lambda_1}2^{(k-1)M}}(z))$
Note that we do not apply the lemma to the larger annulus, where we already have the boundary condition $u=0$ on $\partial B_{\e^{\lambda_1}2^{(\lfloor T/M\rfloor+1)M}}(z)$.

Since $v^k=u$ on $\partial (B_{\e^{\lambda_1}2^{kM}}(z)\setminus B_{\e^{\lambda_1}2^{(k-1)M}}(z))$ we deduce that the function $v$ defined on $\Omega$ by 
$$
v:=
\begin{cases}
v^k & \hbox{ on } B_{\e^{\lambda_1}2^{kM}}(z)\setminus \overline{B}_{\e^{\lambda_1}2^{(k-1)M)}}(z),\, k\in\{1,...,\lfloor T/M\rfloor\} \\
u & \hbox{ otherwise}
\end{cases}
$$
belongs to $H^1_0(\Omega)$ and $v=1$ on $B_\e(z)$. After noting that $ \e^{\lambda_1}2^{(\lfloor T/M\rfloor+1)M}\ge R_\Omega$, we can write
\begin{eqnarray*}
F_{\e,\delta}(v,\Omega)
& =& F_{\e, \delta}(v,B_{\e^{\lambda_1}2^{(\lfloor T/M\rfloor+1)M}}(z)) \\
& =&  F_{\e,\delta}(v,B_{\e^{\lambda_1}}(z))\setminus B_{\e}(z)) + \sum_{k=1}^{\lfloor T/M\rfloor+1}F_{\e,\delta}(v,B_{\e^{\lambda_1}2^{kM}}(z)\setminus B_{\e^{\lambda_1}2^{(k-1)M)}}(z))\\
& =&  F_{\e,\delta}(u,B_{\e^{\lambda_1}}(z))\setminus B_{\e}(z)) + \sum_{k=1}^{\lfloor T/M\rfloor}F_{\e,\delta}(v^k,B_{\e^{\lambda_1}2^{kM}}(z)\setminus B_{\e^{\lambda_1}2^{(k-1)M)}}(z))\\
&& + F_{\e,\delta}(u,B_{\e^{\lambda_1}2^{(\lfloor T/M\rfloor+1)M}}(z)\setminus B_{\e^{\lambda_1}2^{\lfloor T/M\rfloor M)}}(z)),
\end{eqnarray*}
so that, by claim (iv) of Lemma \ref{modlemma1}, we get
\begin{equation}
F_{\e,\delta}(v,\Omega)\le \Bigl(1+\frac{C}{N-1}\Bigr)\, F_{\e,\delta}(u,\Omega).
\end{equation}

Now we make use of the fact that $v$ is a modification of $u$ with the property of being the constant $c_k=c_k^\e$ on spheres centered in $z$ of radius $\e^{\lambda_1}2^{kM-j_k}$, where $j_k\in\{1,...,N-1\}$ for every $k=1,...,\lfloor T/M\rfloor$ (claim (ii) of Lemma \ref{modlemma1}).
We will omit the dependence of $c_k$ on $\e$ as long as this is kept fixed.
Moreover, $v$ preserves the conditions $v=u=1$ on $B_\e(z)$ and $v=u=0$ on $\partial B_{\e^{\lambda_1}2^{(\lfloor T/M\rfloor+1)M}}(z)$. Thus we write
\begin{eqnarray}\label{mainpropch3split}
\nonumber
&&\hskip-2cm F_{\e,\delta}(v,B_{\e^{\lambda_1}2^{(\lfloor T/M\rfloor+1)M}}(z))  = F_{\e, \delta}(v,B_{\e^{\lambda_1}2^{M-j_1}}(z)) \\ \nonumber
&&\qquad\qquad\qquad +\sum_{k=2}^{\lfloor T/M\rfloor}F_{\e,\delta}(v,B_{\e^{\lambda_1}2^{kM-j_k}}(z)\setminus B_{\e^{\lambda_1}2^{(k-1)M-j_{k-1}}}(z)) \\
&&\qquad\qquad\qquad  +\ F_{\e,\delta}(v,B_{\e^{\lambda_1}2^{(\lfloor T/M\rfloor+1)M}}(z)\setminus B_{\e^{\lambda_1}2^{\lfloor T/M\rfloor M-j_{\lfloor T/M\rfloor}}}(z)).
\end{eqnarray}

We now define functions $w^k$ as follows. The function $w^1\in H^1(B_{\e^{\lambda_1}2^M}(z))$ is defined as 
$$
w^1(x)=\begin{cases} v(x) &\hbox{ if $x\in B_{\e^{\lambda_1}2^{M-j_1}}(z)$}\\
c_1 &\hbox{ otherwise.}\end{cases}
$$
In this way 
\begin{equation}\label{est1}
F_{\e,\delta}(w^1,B_{\e^{\lambda_1}2^M}(z))=F_{\e,\delta}(v,B_{\e^{\lambda_1}2^{M-j_1}}(z)).
\end{equation}
For $k\in\{2,...,\lfloor T/M\rfloor+1\}$, we set 
$$
A_{M,k}^N:=B_{\e^{\lambda_1}2^{kM}}(z)\setminus \overline{B}_{\e^{\lambda_1}2^{(k-1)M-N}}(z).
$$
If $k\in\{2,...,\lfloor T/M\rfloor\}$
we define $w^k\in H^1(A_{M,k}^N)$ by setting
$$
w^k(x)=\begin{cases} c_{k-1}&\hbox{ if $x\in{B}_{\e^{\lambda_1}2^{(k-1)M-j_{k-1}}} (z)\setminus \overline{B}_{\e^{\lambda_1}2^{(k-1)M-N}}(z)$}\\
v(x) &\hbox{ if $x\in B_{\e^{\lambda_1}2^{kM-j_k}}(z)\setminus \overline{B}_{\e^{\lambda_1}2^{(k-1)M-j_{k-1}}} (z)$}\\
c_k &\hbox{ if $x\in B_{\e^{\lambda_1}2^{kM}}(z)\setminus \overline{B}_{\e^{\lambda_1}2^{kM-j_k}}(z)$}\end{cases}
$$
For $k=\lfloor T/M \rfloor+1$, we define 
$$
w^{\lfloor T/M \rfloor+1}(x)=\begin{cases} c_{\lfloor T/M \rfloor}&\hbox{ if $x\in{B}_{\e^{\lambda_1}2^{\lfloor T/M \rfloor M-j_{k-1}}} (z)\setminus \overline{B}_{\e^{\lambda_1}2^{(k-1)M-N}}(z)$}\\
v(x) &\hbox{ otherwise.}\end{cases}$$
In this way, we have 
\begin{equation}\label{est2}
F_{\e,\delta}(w^k,A_{M,k}^N)=F_{\e,\delta}(v,B_{\e^{\lambda_1}2^{kM-j_k}}(z)\setminus B_{\e^{\lambda_1}2^{(k-1)M-j_{k-1}}}(z))
\end{equation}
for all $k\in\{2,...,\lfloor T/M\rfloor+1\}$.

If we set $A_{M,1}^N=B_{\e^{\lambda_1}2^M}(z)$, then we can rewrite \eqref{mainpropch3split} as
$$
F_{\e,\delta}(v,B_{\e^{\lambda_1}2^{(\lfloor T/M\rfloor+1)M}}(z)) =  \sum_{k=1}^{\lfloor T/M\rfloor+1} F_{\e,\delta}(w^k,A_{M,k}^N).
$$

By the arbitrariness of  $z\in\Omega$ and $u\in H^1_0(\Omega)$, after setting $c_{\lfloor T/M \rfloor+1}:=0$ it follows
\begin{eqnarray}\nonumber
&&\Bigl(1+\frac{C}{N-1}\Bigr)\,m_{\e,\delta}
\\ \nonumber
&\geq& \Bigl(1+\frac{C}{N-1}\Bigr)\min\Bigl\{F_{\e,\delta}(u,B_{\e^{\lambda_1}2^{(\lfloor T/M\rfloor+1)M}}(z)) :\\ \nonumber
&&\qquad\qquad\qquad u\in H^1_0(B_{\e^{\lambda_1}2^{(\lfloor T/M\rfloor+1)M}}(z)),  u=1 \hbox{ on } B_\e(z)\Bigr\} 
%\end{eqnarray}\begin{eqnarray}
\\ \label{rescalebc}\nonumber
&\geq& \min\Bigl\{F_{\e,\delta}(v,B_{\e^{\lambda_1}2^{M}}(z)) : v\in H^1(B_{\e^{\lambda_1}2^{M}}(z)),\\
&& \qquad\qquad\qquad  \nonumber
v=1 \hbox{ on } B_\e(z), \ v=c_1 \hbox{ on }\partial B_{\e^{\lambda_1}2^{M}}(z)\Bigr\}
\\ \nonumber
&&+\sum_{k=2}^{\lfloor T/M\rfloor+1} \min\Bigl\{F_{\e,\delta}(v, A_{M,k}^N) : v\in H^1(A_{M,k}^N), \\
\nonumber
                   && \qquad\qquad\qquad  v=c_{k-1} \hbox{ on } \partial   B_{\e^{\lambda_1}2^{(k-1)M-N}}(z),  v=c_k \hbox{ on } \partial B_{\e^{\lambda_1}2^{kM}}(z)\Bigr\}
                   \\
   &=&                \min\Bigl\{F_{\e,\delta}(v,B_{\e^{\lambda_1}2^{M}}(z)) : v\in H^1_0(B_{\e^{\lambda_1}2^{M}}(z)), v=1 \hbox{ on } B_\e(z)\Bigr\}\times(1-c_1)^2 
\\ \nonumber
&&+\sum_{k=2}^{\lfloor T/M\rfloor+1} \min\Bigl\{F_{\e,\delta}(v, A_{M,k}^N) : v\in H^1(A_{M,k}^N),  v=1 \hbox{ on } \partial B_{\e^{\lambda_1}2^{(k-1)M-N}}(z), \\ \label{rescalebc-2}
                &&\qquad\qquad v=0 \hbox{ on } \partial B_{\e^{\lambda_1}2^{kM}}(z)\Bigr\}\times(c_{k-1}-c_k)^2, 
\end{eqnarray}
where we have used the 2-homogeneity of minimum problems in the last equality.

The minimum in \eqref{rescalebc} can be simply estimated by 
\begin{eqnarray}\label{eq7}\nonumber
&&\alpha\min\Bigl\{\int_{B_{\e^{\lambda_1}2^{M}}(z)}|\nabla v|^2\,dx : v\in H^1_0(B_{\e^{\lambda_1}2^{M}}(z)), v=1 \hbox{ on } B_\e(z)\Bigr\}
=\alpha\frac{2\pi}{\log (\e^{(\lambda_1-1)}2^M)}\\ &=&\frac{2\pi\alpha}{(1-\lambda_1)|\log\e|+M \log 2}\,.
 \end{eqnarray}
As for the minima in \eqref{rescalebc-2}, with fixed $k\in\{2,...,\lfloor T/M\rfloor+1\}$ the change of variable $ w(x)= v(z+x\,\e^{\lambda_1}2^{(k-1)M})$ leads to rewriting the corresponding minimum as
\begin{eqnarray*}\label{eq8}\nonumber
&&\min\Bigl\{F_{\e,\delta}(v, A_{M,k}^N) : v\in H^1(A_{M,k}^N), \\ \nonumber&&\qquad\qquad\qquad\qquad 
v=1 \hbox{ on } \partial B_{\e^{\lambda_1}2^{(k-1)M-N}}(z)),
                 v=0 \hbox{ on } \partial B_{\e^{\lambda_1}2^{kM}}(z)\Bigr\}
\\ \nonumber
&=&\min\Bigl\{\int_{B_{2^M}\setminus B_{2^{-N}}}  a\left(\frac{x}{\delta}\e^{\lambda_1}2^{(k-1)M}+\tau_k^\e\right)|\nabla w(x)|^2\,dx : w\in H^1(B_{2^M}\setminus \overline{B}_{2^{-N}}), \\
&&\qquad\qquad\qquad\qquad w=1 \hbox{ on } \partial B_{2^{-N}}, w=0 \hbox{ on } \partial B_{2^M}\Bigr\},
 \end{eqnarray*}
 where $\tau_k^\e= -\frac{z}{\delta}\e^{\lambda_1}2^{(k-1)M}$. We can use Theorem \ref{homab} with $A=B_{2^M}\setminus \overline B_{2^{-N}}$,
 $\eta= {\delta\over\e^{\lambda_1}2^{(k-1)M}}$, noting that $\eta\to0$ by \eqref{eqz}, and $\phi$ any function in $H^1(B_{2^M}\setminus \overline{B}_{2^{-N}})$ such that $\phi=1$ on $\partial B_{2^{-N}}$ and $\phi=0$ on $\partial B_{2^M}$, so that this last minimum is estimated by  
\begin{eqnarray*}
&&\hskip-1cm \sqrt{\alpha\beta} \min\Bigl\{\int_{B_{2^M}\setminus B_{2^{-N}}}|\nabla w(x)|^2\,dx : w\in H^1(B_{2^M}\setminus \overline{B}_{2^{-N}}): \\ &&
\qquad\qquad\qquad\qquad w=1 \hbox{ on } \partial B_{2^{-N}}, v=0 \hbox{ on } \partial B_{2^M}\Bigr\}\\
&=& 2\pi \sqrt{\alpha\beta}{1\over (M+N)\log2}
\end{eqnarray*}
up to an infinitesimal term as $\e\to0$. This term can be chosen independent of $k$ since $\eta\le {\delta\over\e^{\lambda_1}}2^{M}$ for all $k$ and ${\delta\over\e^{\lambda_1}}2^{M}$ is infinitesimal as $\e\to0$.

By the convexity inequality for the square and the fact that $\sum_{k=2}^{\lfloor T/M \rfloor+1}(c_{k-1}-c_k)= c_1$,  we have 
$$
\sum_{k=2}^{\lfloor T/M \rfloor+1}(c_{k-1}-c_k)^2 \ge \frac{c_1^2}{\lfloor T/M\rfloor-1}\ge \frac{Mc_1^2}{T}\ge \frac{c_1^2 M \log 2}{\lambda_1 |\log\e|+\log2},
$$
having taken \eqref{T} into account.

Gathering the estimates above, we obtain
\begin{eqnarray*}
&&\Bigl(1+\frac{C}{N-1}\Bigr)m_{\e,\delta}
\\
%&\ge& \frac{2\pi \alpha(1-c_1)^2}{(1-\lambda_1)|\log \e|+M \log 2}+\Bigl(\frac{2 \pi \sqrt{\alpha\beta}}{(M+N)\log2}+o_{\e}(1)\Bigr)
%\sum_{k=2}^{\lfloor T/M \rfloor+1}(c_{k-1}-c_k)^2 
%\\
%&\ge& \frac{2\pi \alpha(1-c_1)^2}{(1-\lambda_1)|\log \e|+M \log 2}+\Bigl(\frac{2 \pi \sqrt{\alpha\beta}}{(M+N)\log2}+o_{\e}(1)\Bigr)
% \frac{c_1^2}{\lfloor T/M\rfloor-1}\\
%&\ge& \frac{2\pi \alpha(1-c_1)^2}{(1-\lambda_1)|\log \e|+M \log 2}+\Bigl(\frac{2 \pi \sqrt{\alpha\beta}}{(M+N)\log2}+o_{\e}(1)\Bigr)
% \frac{Mc_1^2}{T}
% \\
&\ge& \frac{2\pi \alpha(1-c_1)^2}{(1-\lambda_1)|\log \e|+M \log 2}+\Bigl(\frac{2 \pi \sqrt{\alpha\beta}}{(M+N)\log2}+o_{\e}(1)\Bigr)
 \frac{c_1^2 M \log 2}{\lambda_1 |\log\e|+\log2}.
\end{eqnarray*}
Reinstating the dependence of $c_1$ on $\e$ in the notation, we deduce that 
\begin{eqnarray}\label{finalest}\nonumber
\Bigl(1+\frac{C}{N-1}\Bigr)|\log\e|m_{\e,\delta} & \geq& \frac{2\pi \alpha(1-c_1^\e)^2|\log\e|}{(1-\lambda_1)|\log \e|+M \log 2} \\
&& + \biggl[\frac{2 \pi \sqrt{\alpha\beta}}{(M+N)\log2}+o_{\e}(1)\biggr]\frac{(c_1^\e)^2 M \log 2 |\log\e|}{\lambda_1 |\log\e|+\log2}\,,
\end{eqnarray}
for every $N<M<T$. We can assume that $c_1^\e\to c_1$ as $\e\to 0$. Since $T\to+\infty$ as $\e\to0$, we get
$$
\Bigl(1+\frac{C}{N-1}\Bigr)\liminf_{\e\to0}|\log\e|m_{\e,\delta} \geq \frac{2\pi \alpha(1-c_1)^2}{(1-\lambda_1)} + \frac{2 \pi \sqrt{\alpha\beta}\,c_1^2}{\lambda_1}\frac{M}{M+N}\,,
$$
for every $N,M$ natural numbers such that $N<M$.
We successively let $M\to+\infty$ and then $N\to+\infty$ obtaining
\[
\liminf_{\e\to0}|\log\e|m_{\e,\delta}\geq \frac{2\pi \alpha(1-c_1)^2}{(1-\lambda_1)} + \frac{2 \pi \sqrt{\alpha\beta}\,c_1^2}{\lambda_1}.
\]

This expression in $c_1$ is minimized for
\begin{equation}\label{cilambda1}
c_1=\frac{\alpha\lambda_1}{\alpha\lambda_1+ \sqrt{\alpha\beta}(1-\lambda_1)},
\end{equation}
which gives the minimal value
\begin{equation}\label{cilambdamin}
2\pi\frac{\alpha\sqrt{\alpha\beta}}{\lambda_1\alpha+(1-\lambda_1)\sqrt{\alpha\beta}}\,.
\end{equation}
We finally observe that our estimate holds true for every $\lambda_1<\lambda$, thus 
\begin{equation}\label{lb}
\liminf_{\e\to0}\,|\log\e|m_{\e,\delta} \geq \sup_{0<\lambda_1<\lambda}2\pi\frac{\alpha\sqrt{\alpha\beta}}{\lambda_1\alpha+(1-\lambda_1)\sqrt{\alpha\beta}}=2\pi\frac{\alpha\sqrt{\alpha\beta}}{\lambda\alpha+(1-\lambda)\sqrt{\alpha\beta}},
\end{equation}
which is the desired lower bound.

\subsection{Construction of optimal capacitary profiles}
We now construct functions (almost) realizing the lower bound \eqref{lb}. The corresponding small balls will be centered in a point where the value of $a$ is $\alpha$ uniformly distant from the boundary of $\Omega$. To that end,
for every $\e>0$ we consider $z_\e\in\Omega$ centre of a square $Q_{\e}$ of side length $\delta/2$ in which $a(x/\delta)=\alpha$ and such that there exists $R_0<1$ independent of $\e$ such that $B_{R_0}(z_\e)\subset\Omega$. We now subdivide the construction in the case $\lambda=0$, in which it is trivial, and $\lambda\in(0,1]$.

\bigskip
If $\lambda=0$ then we choose 
$$
u_{\e}(x)=
\begin{cases}
1 & \hbox{ if } |x-z_{\e}| \leq \e \\
1-{1\over \log\delta-\log\e}\log\frac{|x-z_{\e}|}{\e} & \hbox{ if }\e < |x-z_{\e}| <\delta\\
0 &  \hbox{ otherwise.}
\end{cases}
$$
A direct computation shows that 
$$
F_{\e,\delta}(u_\e,\Omega)=2\pi \alpha{1\over |\log\e-\log \delta|},
$$
which implies that
$$
\limsup_{\e\to0}|\log\e|m_{\e,\delta} \leq 2\pi\alpha,
$$
and the claim.

\bigskip

Let $\lambda\in (0,1]$ and let $$
c:=c(\lambda)=\frac{\alpha\lambda}{\alpha\lambda+{\sqrt{\alpha\beta}}(1-{\lambda})}
$$
If $\lambda\in (0,1)$ this is the optimal value in \eqref{cilambda1}, while for $\lambda=1$ we have $c=1$.
We consider $\lambda_1$ with $0<\lambda_1< \lambda$, so that $\delta<\!< \e^{\lambda_1}$ as $\e\to0$, and let $u_{\e}^0\in H^1(B_{\e^{\lambda_1}}(z_{\e}))$ be defined by
$$
u_{\e}^0(x)=
\begin{cases}
1 & \hbox{ if } |x-z_{\e}| \leq \e \\
1-\frac{1-c}{(1-\lambda_1)|\log\e|}\log\frac{|x-z_{\e}|}{\e} & \hbox{ if }\e < |x-z_{\e}| \leq \e^{\lambda_1}. 
\end{cases}
$$
Note that $u_{\e}^0$ is the solution of 
$$
\min\Bigl\{\int_{B_{\e^{\lambda_1}}(z_{\e})}|\nabla u(x)|^2\,dx : u\in H^1(B_{\e^{\lambda_1}}(z_{\e})),  u=1 \hbox{ on } B_\e(z_{\e}), u=c \hbox{ on } \partial B_{\e^{\lambda_1}}(z_{\e}) \Bigr\}.
$$

\smallskip
We now define an optimal capacitary profile outside $B_{\e^{\lambda_1}}(z_{\e})$ using a construction on concentric annuli  homothetic to $B_1\setminus \overline B_{1/2}$. To that end, for all $\eta>0$, $\tau_\eta$ and $F_\eta$ as in Theorem \ref{homab} with $A= B_1\setminus\overline  B_{1/2}$, we consider the solution $w_\eta$ of the problem 
\begin{equation}
m_\delta=\min\Bigl\{F_\eta(w): w=1 \hbox{ on } \partial B_{1/2},  w=0 \hbox{ on } \partial B_1\}.
\end{equation}
Note that by Theorem \ref{homab} there exists $\omega:[0,+\infty)\to[0,+\infty)$ increasing and with $\omega(0)=0$ such that
\begin{equation}
\Bigl|m_\delta-{2\pi\sqrt{\alpha\beta}\over\log 2}\Bigr|\le \omega(\eta).
\end{equation}

We now let $T:=\max\{n\in \mathbb{N} : \e^{\lambda_1}2^n \leq R_0\}$. For $\e>0$ small enough $T$ is well defined and $$T=\Bigl\lfloor\frac{\lambda_1|\log\e|+\log R_0}{\log2}\Bigr\rfloor,$$ so that we can assume $T\geq2$. 
We set
$$A_k=B_{\e^{\lambda_1}2^k}(z_{\e})\setminus \overline{B}_{\e^{\lambda_1}2^{k-1}}(z_{\e})$$
for all $1\le k\le T$.

For $k$ fixed, let $\eta= {\delta\over e^{\lambda_1}2^{k}}$, let $w_\eta$ be defined above and let $u_{\e}^k\in H^1(A_k)$ be defined by 
$$
u_{\e}^k(x)= \frac{T-k}{T}\,c+ {c\over T} w_\eta ((x-z)\e^{-\lambda_1}2^{-k})\,,
$$
so that  $u_{\e}^k$ is the solution of the minimum problem
\begin{eqnarray*}
m^k_\e&=&\min\Bigl\{F_{\e,\delta}(v, B_k) : \, v\in H^1(A_k),  v=\hbox{${T-k+1\over T}$}\,c \hbox{ on } \partial B_{\e^{\lambda_1}2^{k-1}}(z_{\e}), \\
&&\qquad v=\hbox{$\frac{T-k}{T}$}\,c \hbox{ on } \partial B_{\e^{\lambda_1}2^k}(z_{\e} ) \Bigr\} \\
&= &{c^2\over T^2}m_\eta ,
\end{eqnarray*}
so that
$$
\Bigl|m^k_\eta-{c^2\over T^2}{2\pi\sqrt{\alpha\beta}\over\log 2}\Bigr|\le {c^2\over T^2}\omega(\eta)\le {c^2\over T^2}\omega\Bigl({\delta\over \e^{\lambda_1}}\Bigr)
$$

Define
\[
u_{\e}:=
\begin{cases}
u_{\e}^0 & \hbox{ on } B_{\e^{\lambda_1}}(z_{\e}) \\
u_{\e}^k & \hbox{ on } A_k,\, k\in\{1,...,T\} \\
0 & \hbox{ on } \Omega\setminus B_{\e^{\lambda_1}2^T}(z_{\e}),
\end{cases}
\]
which is an admissible function for the problem defining $m_{\e,\delta}$.

We now proceed with the estimates of $F_{\e,\delta}(u_\e,\Omega)$. If $\lambda=1$ then $c=1$ and $F_{\e,\delta}(u_{\e}, B_{\e^{\lambda_1}}(z_{\e}))=0$. If $\lambda\in (0,1)$ let 
 $\lambda_2$ with $\lambda< \lambda_2<1$.
 This condition ensure that $\e^{\lambda_2} <\!< \delta<\!< \e^{\lambda_1}$ as $\e\to0$, so that we deduce the inclusions $$
B_{\e^{\lambda_2}}(z_{\e})\subseteq Q_{\e}\subseteq B_{\e^{\lambda_1}}(z_{\e}).
$$
Hence, we have
\begin{equation}\label{minseq1}
\begin{split}
F_{\e,\delta}(u_{\e}, B_{\e^{\lambda_1}}(z_{\e})) & = 
\int_{B_{\e^{\lambda_2}}(z_{\e})}a\Bigl(\frac{x}{\delta}\Bigr)|\nabla u^0_{\e}|^2\,dx + 
\int_{B_{\e^{\lambda_1}}(z_{\e})\setminus B_{\e^{\lambda_2}}(z_{\e})} a\Bigl(\frac{x}{\delta}\Bigr)|\nabla u^0_{\e}|^2\,dx \\
& \leq \alpha\int_{B_{\e^{\lambda_2}}(z_{\e})}|\nabla u^0_{\e}|^2\,dx + \beta\int_{B_{\e^{\lambda_2}}(z_{\e}))\setminus B(z_{\e},\e^{\lambda_2})}|\nabla u^0_{\e}|^2\,dx \\
%& = 2\pi\biggl(\frac{1-c}{(1-\lambda_1)|\log\e|}\biggr)^2\biggl(\alpha \int_{\e}^{\e^{\lambda_2}}\frac{1}{\rho}\,d\rho + \beta\int_{\e^{\lambda_2}}^{\e^{\lambda_1}}\frac{1}{\rho}\,d\rho \biggr) \\
%& = 2\pi\biggl(\frac{1-c}{(1-\lambda_1)|\log\e|}\biggr)^2\biggl(\alpha(1-\lambda_2)|\log\e|+\beta(\lambda_2-\lambda_1)|\log\e|\biggr) \\
& = \frac{2\pi(1-c)^2}{|\log\e|}\Bigl(\alpha\frac{1-\lambda_2}{(1-\lambda_1)^2}+\beta\frac{\lambda_2-\lambda_1}
{(1-\lambda_1)^2} \Bigr).
\end{split}    
\end{equation}

As for the remaining part, %since $T\geq\lfloor\frac{\lambda_1|\log\e|+\log dist(\partial B,\partial\Omega)}{\log2}\rfloor\geq \frac{\lambda_1|\log\e|+\log dist(\partial B,\partial\Omega)}{\log2}-1$, we deduce
we can write it as
\begin{equation}\label{minseq2}
\begin{split}
    \sum_{k=1}^T F_{\e,\delta}(u_{\e},B_k) & = \Bigl(\sqrt{\alpha\beta}\frac{2\pi}{\log 2}+o_{\e}(1) \Bigr)\frac{c^2}{T} \\
    & \leq \biggl(\sqrt{\alpha\beta}\frac{2\pi}{\log 2}+o_{\e}(1) \biggr) \frac{c^2\log2}{\lambda_1|\log\e|+\log R_0}.
\end{split}
\end{equation}

Combining $\eqref{minseq1}$ and $\eqref{minseq2}$ we get
\[
\begin{split}
|\log\e|m_{\e,\delta} & \leq 2\pi(1-c)^2\biggl(\alpha\frac{1-\lambda_2}{(1-\lambda_1)^2}+\beta\frac{\lambda_2-\lambda_1}{(1-\lambda_1)^2}\biggr) \\
& +\biggl(\sqrt{\alpha\beta}\frac{2\pi}{\log 2}+o_{\e}(1) \biggr) \frac{c^2\log2|\log\e|}{\lambda_1|\log\e|+\log R_0}\,.
\end{split}
\]
Letting $\e\to 0$ we then obtain
$$
\limsup_{\e\to0}|\log\e|m_{\e,\delta}\le 2\pi(1-c)^2\biggl(\alpha\frac{1-\lambda_2}{(1-\lambda_1)^2}+\beta\frac{\lambda_2-\lambda_1}{(1-\lambda_1)^2}\biggr)+2\pi \sqrt{\alpha\beta} \frac{c^2}{\lambda_1}
$$

If $\lambda=1$ then $c=1$, so that the first term in the right-hand side term is $0$ and, after letting $\lambda_1\to1$, we have 
$$
\limsup_{\e\to0}|\log\e|m_{\e,\delta}\le 2\pi \sqrt{\alpha\beta}.$$

If $\lambda\in (0,1)$ then we let $\lambda_1,\lambda_2\to \lambda$, so that we get 
$$
\limsup_{\e\to0}|\log\e|m_{\e,\delta}\le 2\pi(1-c)^2\alpha\frac{1}{1-\lambda}+2\pi \sqrt{\alpha\beta} \frac{c^2}{\lambda}\,.
$$
Recalling \eqref{cilambdamin}, we get 
$$
\limsup_{\e\to0}|\log\e|m_{\e,\delta} \leq 2\pi\frac{\alpha\sqrt{\alpha\beta}}{\lambda\alpha+(1-\lambda)\sqrt{\alpha\beta}},
$$
and the claim.

\subsection{Remarks}
We briefly comment the case of arbitrary measurable $a$. If $a$ is continuous then formula \eqref{casoii-0} rereads
\begin{equation}\label{casoii-1}
\lim_{\e\to 0}|\log\e| m_{\e,\delta}=%\begin{cases} \alpha &\hbox{ if } \lambda=0\\
2\pi {\alpha'\beta'\over \lambda \alpha'+(1-\lambda)\beta'}, 
%&\hbox{ if } 0<\lambda<1\\C_{\rm hom} &\hbox{ if } \lambda=1\end{cases}
\end{equation}
where $\alpha'=\min a$ and $\beta'=\sqrt{{\rm det}A_{\hom}}$, $\langle A_{\hom}\xi,\xi\rangle$ being the integrand of the $\Gamma$-limit as in Theorem \ref{homab}. The proof is the same, with small changes in the computation, taking balls centered in $z_\e$ minimizing $a$ for the construction and using Theorem \ref{homab}. A minor issue is the fact that if 
$A_{\hom}$ is not a multiple of the identity the computation in \eqref{eqRr} is not exact but must be slightly adapted (e.g.~using ellipses in the place of circles).

If $a$ is only measurable then the formula holds with $\alpha'$ the essential minimum of $a$, and some extra care must be used as in the analog result in \cite{ABCDP}. We may recover the checkerboard case above noting that $A_{\hom}$ is equal to $\sqrt{\alpha\beta}$ times the identity, so that det$A_{\hom}=\alpha\beta$.

For an extension to any dimension or to general integrands, we cannot use the special structure of quadratic forms and some issues are more technical. We refer to \cite{Brusca} for an asymptotic formula and applications in that case.

\subsection*{Acknowledgments}
The problem has been inspired by a conversation with Fran\c cois Murat.
AB is a member of GNAMPA, INdAM.

\end{document}